\documentclass[11pt]{article}

\usepackage{amsfonts,amsthm,amssymb}
\usepackage[cp866]{inputenc}

\usepackage[T2A]{fontenc}

\usepackage{amsmath,amsthm,amssymb}
\usepackage{mathtext}

\usepackage{epsfig}
\usepackage{graphicx}
\usepackage{wrapfig}

\begin{document}

\newtheorem{lm}{Lemma}
\newtheorem{theorem}{Theorem}
\newtheorem{prop}{Proposition}
\newtheorem{df}{Definition}
\newtheorem{remark}{Remark}
\newtheorem{corollary}{Corollary}
\newtheorem{ex}{Example}

\begin{center}{\large\bf On Shilnikov attractors of three-dimensional flows and maps.}
\end{center}
\medskip
\begin{center}{\bf Bakhanova Yu.V.$^1$, Gonchenko S.V.$^{1,2}$, Gonchenko A.S.$^2$, \\ Kazakov A.O.$^{1,2}$, Samylina E.A.$^1$ }
\end{center}
\medskip
\begin{center}{$^2$ National Research University Higher School of Economics, Nizhny Novgorod, Russia}
\end{center}
\begin{center}{$^1$ Lobachevsky State University of Nizhny Novgorod, Russia}
\end{center}

\begin{center} {e-mail: bakhanovayu@gmail.com, phone: +79108729548; e-mail: sergey.gonchenko@mail.ru, phone:+79065567766; e-mail: agonchenko@mail.ru, phone:+79527617397;  e-mail: kazakovdz@yandex.ru, phone: +79081659694; e-mail: samylina\_evgeniya@mail.ru, phone: +79527734356}
\end{center}
\bigskip

\begin{abstract}
We describe scenarios for the emergence of Shilnikov attractors, i.e. strange attractors containing a saddle-focus with two-dimensional unstable manifold,  in the case of three-dimensional flows and maps. The presented results are illustrated with various specific examples.
\end{abstract}


\section{Introduction}

The discovery of spiral chaos, i.e. the proof of the existence of a complex structure of orbits in a neighbourhood of a homoclinic loop of a saddle-focus equilibrium state, is rightfully considered one of the most significant achievements of the theory of dynamical systems. This discovery was made by L.P. Shilnikov in his work \cite{Sh65}, published in 1965.
At that time, in both mathematics and physics, there was practically no more or less suitable concept for explaining such phenomena in models described by finite-dimensional deterministic systems. Therefore, the complex structure (chaos) in neighbourhoods of homoclinic loops of  equilibria was absolutely unexpected and completely contradicted the concepts based on the theory of two-dimensional systems. Moreover, such significant difference in the orbit structure
for a saddle and a saddle-focus was seemed very strange, since from the topological point of view, the saddle and saddle-focus are arranged in the same way.
However, they differ globally. Thus, from the Shilnikov conditions of nondegeneracy for a homoclinic loop of a saddle \cite{Sh62,Sh63,Sh68},
one can deduce that in its neighbourhood there is a smooth two-dimensional global central invariant manifold, and thus the problem is effectively
two-dimensional. Whereas the saddle-focus does not have such a manifold even locally, and therefore the dynamics here is principally multidimensional.

A saddle-focus differs from a saddle in that, among its leading eigenvalues (those nearest to the imaginary axis), a pair of complex-conjugate ones exists. Then in the case of three-dimensional systems saddle-focus equilibria are divided into two different types:
\begin{itemize}
\item
a {\em saddle-focus (2,1)} with two-dimensional stable and one-dimensional unstable manifolds that has eigenvalues $\lambda \pm i \omega, \gamma $, where $\lambda<0, \gamma>0, \omega \neq 0$;
\item
a {saddle-focus (1,2)} with one-dimensional stable and two-dimensional unstable manifolds that has eigenvalues $\lambda, \gamma \pm i \omega $, where $\lambda<0, \gamma>0, \omega \neq 0$.
\end{itemize}
Back in the papers \cite{Sh62, Sh63} Shilnikov showed that bifurcations of  systems with a homoclinic loop of a saddle-focus (2,1) with negative saddle value
$$
\sigma = \lambda + \gamma
$$
do not differ from the case of a saddle -- here a unique stable limit cycle is born when the loop spits inward  (this is due to the fact that the first return map turns out to be contracting in this case). The case of a three-dimensional saddle-focus (1,2) is obviously reduced to the case of a saddle-focus  (2,1) by time reversal.
Therefore, here in the case  $\sigma>0$  a unique completely unstable limit cycle is born from a homoclinic loop of a saddle-focus  of type (1,2).

However, if the closest to the imaginary axis eigenvalues are complex conjugate ones, i.e., if $\sigma> 0 $ in the case of a saddle-focus (2,1) or $\sigma <0$ in the case of a saddle-focus (1,2), the situation becomes completely different -- already the system itself with a homoclinic loop has a complex structure, since, in any neighbourhood of the loop, there are infinitely many saddle periodic orbits \cite{Sh65}.
In fact, Shilnikov discovered in \cite{Sh65} that the existence of a homoclinic loop of such saddle-focus implies chaos. The notion itself
did not exist then (in 1965); the ``chaos theory'' emerged and became popular only 10-20 years later. Chaos was
found in many nonlinear models
and it also
occurred that strange attractors in models of various origins often have a spiral structure, i.e. the
chaotic orbits seem to move near a saddle-focus homoclinic loop.

According to the classification of saddle-foci, we will also distinguish two types of spiral homoclinic attractors of three-dimensional flows:

\begin{itemize}
\item[1)] a {\em figure-8 spiral attractor}, when the attractor contains a saddle-focus (2,1) (with $\sigma>0$) and entirely both its one-dimensional unstable separatrices (composing a  homoclinic-8);
\item[2)] a {\em Shilnikov attractor}, when the attractor contains a saddle-focus (1,2) (with $\sigma<0$) and entirely its two-dimensional unstable manifold.
\end{itemize}

That the homoclinic loop to a saddle-focus implies chaos -- this is the Shilnikov theorem \cite{Sh65,Sh70}, but why is the converse also so often true, i.e. why is the observed chaos often spiral? This question was of great interest to Shilnikov.
He discovered \cite{Sh84,Sh86} that if a system depends on a parameter and evolves as it changes, from a stable (stationary) regime to a chaotic one, then on this way a saddle-focus equilibrium arises naturally and, moreover, its stable and unstable manifolds can come close enough to each other, so that the creation of a homoclinic loop (and chaos) becomes quite expected.

Shilnikov described the corresponding scenarios in the paper \cite{Sh86} for the case of a one-parameter family $X_\mu : \dot x = X(x,\mu)$ of multidimensional systems.
In Section~\ref{sec:Shiscen_flow} we discuss the Shilnikov scenario for the three-dimensional case. As one can see, this scenario
is extremely simple (if to ignore certain fine intermediate details that themselves can be very complicated), and therefore, it should come as no surprise that it is often seen in many models.\footnote{Their list is quite large, we will indicate only some of the most known three-dimensional models.
Thus, spiral chaos was found in radio-electronic devices, such as the Chua circuit \cite{ChuaKomMat1986}, the Anishchenko-Astakhov generator \cite{Ani1990}, in the optical laser systems \cite{ArrMeuGad1987,ArrLapMeuRov1988,PisurtMeuBrechou3,ZhouKurtAllBocMeu2003}, in chemical systems \cite{ArgArnRich1987,ArgArnEleRich1993}, in a certain class of models describing the behavior of neurons \cite{FeudelPei2000}, in biophysical experiments \cite{PartEdGri2003}, in electromechanical systems \cite{KopGasSlu1992n2001,CheWoaDomn2001}, in electrochemical processes \cite{BasHud1988, Noh2009}, in nonlinear convection in magnetic fields \cite{Rucklidg1993}, in mechanical systems \cite{HenLeviOdeh1991}, etc.}
We consider some of these models in Section~\ref{sec:exflows} as examples of systems in which Shilnikov's scenario is implemented.

It is interesting to note that spiral chaos was not initially detected in a number of such models, the studies of chaotic dynamics ended  either at the stage of Feigenbaum-type attractors (which are formed inside the Shilnikov whirlpool (see Sec.~\ref{sec:Shiscen_flow}) as a result of an infinite cascade of period-doubling bifurcations with limit cycles) or a little further, when the attractor in the Poincar\'e section becomes  ``one-component'', like the R\"ossler attractor (an example of such attractor is shown in Fig.~\ref{fig:Arneodo_scenario}e), or the interpretation of the process of its occurrence was inadequate etc.  Of course, this situation was a consequence of the lack of a ``guiding thread'' -- the Shilnikov scenario in this case.

In the work \cite{Sh86} Shilnikov expressed one more important idea that the scenario he proposed for flows can be generalized to the case of maps. In this case, the generalization is direct: the equilibrium state is replaced by a fixed (or periodic) point, and the local and global bifurcations involved in the scenario are replaced by their analogs  for maps. In Section~\ref{sec:Shiscen_map} we give a phenomenological description of such scenarios for the case of orientable and nonorientable three-dimensional maps, and in Section~\ref{sec:exShattr} we give some examples of how such scenarios can be implemented in the case of  three-dimensional generalized H\'enon maps.

\section{On Shilnikov scenario for three-dimensional flows} \label{sec:Shiscen_flow}

\begin{figure} [htb]
	\centerline{
		\includegraphics[width=14cm]{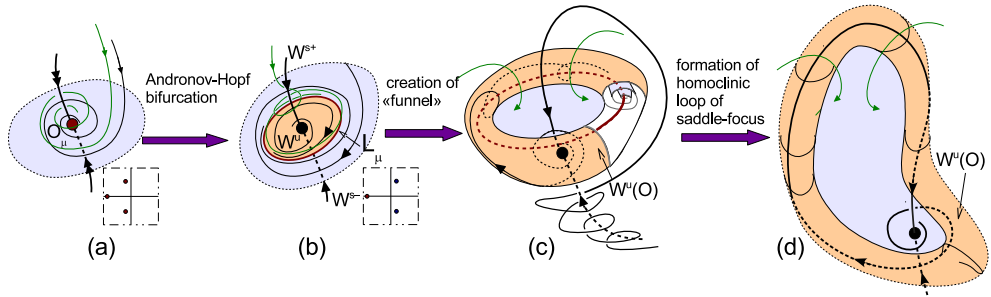}}
	\caption{{\footnotesize The main stages of the phenomenological scenario for the emergence of the Shilnikov attractor. One of the criteria for the strangeness of such an attractor is the existence of a homoclinic (double-asymptotic) orbit to  the saddle-focus $O_\mu$.}}
	\label{figShscen}
\end{figure}

Consider a one-parameter family $X_\mu : \dot x = X(x,\mu)$ of three-dimensional flows such that the system $X_\mu$ has at $\mu_0<\mu<\mu_1$ a single stable equilibrium $O_\mu$ in some absorbing domain $D$, see Fig.~\ref{figShscen}a. In other words, the system $X_\mu$ for $\mu_0<\mu<\mu_1$ has in $D$ a global attractor that is the stable equilibrium $O_\mu$. Let at $\mu \in (\mu_0,\mu_1)$ and close to $\mu_1$ the equilibrium $O_\mu$ have eigenvalues $\nu_{1,2} = \lambda(\mu) \pm i \omega(\mu), \nu_3 =\lambda_{ss}(\mu)$, where $\lambda_{ss} < \lambda <0 $. This means that $O_\mu$ is a stable focus by its leading eigenvalues $\nu_1$ and $\nu_2$.

Assume that at $\mu = \mu_1$ a \textit{supercritical} (soft) Andronov-Hopf bifurcation occurs with $O_\mu$. Then the eigenvalues $\nu_{1,2}$ evolve as $\mu$ changes in such a way that $\lambda(\mu) <0$ if $\mu < \mu_1$, $\lambda(\mu_1) =0$,  and $\lambda(\mu) >0$ if $\mu > \mu_1$. Thus, at $\mu> \mu_1$ the equilibrium $O_\mu $ becomes a saddle-focus (1,2) and a stable limit cycle $L_\mu $ is born in its neighbourhood. Accordingly, for all sufficiently small $\mu >  \mu_1$, the global attractor in $D$ of the system $X_\mu$ is this asymptotically stable limit cycle $L_\mu$. The type of stability of the cycle
$L_\mu$ is determined by its multipliers $\rho_1(\mu)$ and $\rho_2 (\mu)$, which are initially both positive and less than one. In this case, the unstable manifold $W^u(O_\mu)$ of the point $O_\mu$ is a two-dimensional disk with the boundary $L_\mu$, and orbits starting from points on $W^u(O_\mu) $ are wound to $L_\mu$ along spirals ``from the inside'', see Fig.~\ref{figShscen}b.

We assume that with a further increase in $\mu$ the multipliers
$\rho_1(\mu)$ and $\rho_2(\mu)$ of $L_\mu$ become equal at some $\mu = \mu^*$  and complex conjugate for $\mu> \mu^*$. This ``smooth bifurcation'' corresponds to change of type of $L_\mu$, from nodal type for $\mu< \mu^*$ to focal one for $\mu> \mu^*$. Then, the two-dimensional unstable manifold $W^u(O_\mu)$ of the equilibrium $O_\mu$ begins to twist onto the cycle $L_\mu$. In this case,  $W^u(O_\mu)$ takes a form of roulette which is a boundary of the so-called ``Shilnikov whirlpool'', inside of which all orbits of system $X_\mu$ starting in $D$ are drawn (generally speaking, except for one orbit -- the stable separatrix $W^{s-}(O_\mu)$ of the equilibrium $O_\mu$), see
Fig.~\ref{figShscen}c.

When $\mu$ changes,  the sizes of the whirlpool increase, the ``retraction inward of orbits'' is preserved, but the limit cycle $L_\mu$ may lose its stability. In particular, a strange attractor can form in its place as a result of a series of various bifurcations\footnote{The sequence of these bifurcations can be very diverse: this can be a cascade of period doubling bifurcations followed by the appearance of  a Feigenbaum attractor or, few further,  a H\'enon-like attractor in the Poincar\'e map; or a two-dimensional stable invariant torus
$T_\mu$ can be born from the cycle, which can then break down, for example, according to the Afraimovich-Shilnikov scenario \cite{ASh83}, giving rise a strange attractor of the ``torus-chaos'' type; some limit cycle can undergo a subcritical Andronov-Hopf or period-doubling bifurcation after which a chaotic dynamics inside the whirlpool appears instantly; etc. }. This attractor, in turn, can transform into the Shilnikov attractor, i.e. a strange attractor that contains the saddle-focus $O_\mu$ of type (1,2) and entirely its two-dimensional unstable manifold.

The moment $\mu = \mu_{\ell}$ of the formation of a homoclinic loop of saddle-focus $O_\mu$ can be traced by observing for the stable separatrix $W^{s+}(O_\mu)$: at  $\mu < \mu_{\ell}$ it leaves $D$ (as in Fig.~\ref{figShscen}c), at  $\mu > \mu_{\ell}$ it enters the whirlpool, and, at $\mu = \mu_{\ell}$, it falls on the unstable manifold $W^u(O_\mu)$ (thus, $W^{s+}(O_\mu)$ entirely belongs to $W^u(O_\mu)$).

With a further change in $\mu$, the separatrix $W^{s+}(O_\mu)$ entirely enters the whirlpool, makes many turns there, then can again lie on $W^u(O_\mu)$ forming a multi-round homoclinic loop and etc. The moments of formation of such loops are discrete in the parameter $\mu$, but they are not isolated: for each such value of the parameter $\mu$, those values of $\mu$ are accumulated again that correspond to secondary multi-round loops \cite{GGNT97}.

Discrete moments with respect to the parameter $\mu $, when the separatrix $W^{s+}(O_\mu)$ forms homoclinic loops, correspond to the existence of Shilnikov homoclinic attractors\footnote{ of course, when inside the whirlpool there are no other attractors -- for example, local bifurcations may lag behind global bifurcations  and then,  it can happen that both the limit cycle $L_\mu$ is stable and a homoclinic loop exists}. At these moments, the equilibrium $O_\mu$ enters the attractor together with its unstable manifold $W^u(O_\mu)$, and typical orbits of the attractor visit any arbitrarily small neighborhood of the point $O_\mu$. The latter property makes it relatively easy to find the moments of formation of Shilnikov homoclinic attractors using, for example, automated plotting of the $\mu$-dependence of the distance of attractor points from the saddle-focus.

\section{Examples of Shilnikov attractors in three-dimensional flows} \label{sec:exflows}

For a long time the Shilnikov's works \cite{Sh65,Sh70,Sh67}  on  spiral chaos were practically unknown to the mathematical community. They gained a world recognition due to a series of quite affordable papers by Arneodo, Coullet, Tresser, and Spiegel \cite{ArnCoulTres1980,ArnCoulTres1981,ArnCoulTres1982,ArnCoulTres1985}, in which the importance of the Shilnikov homoclinic loop for the chaos theory was emphasized. In particular,
in \cite{ArnCoulTres1981}, geometric illustrations of the Shilnikov's theorem were given. Also, in this paper, it was presented a simple example of a piece-wise linear oscillatory system, in which the existence of a homoclinic loop of saddle-focus (2,1) was analytically established and figure-8 spiral strange attractors were found numerically.

In this section we consider two examples of three-dimensional system in which the Shilnikov attractors are observed: the first example is the well-known ACT-system from the paper \cite{ArnCoulTres1982}, and the second example is the Gaspard--Nicolis model of chemical oscillator from \cite{GaspardNicolis83}.

\subsection{Shilnikov attractor in an Arneodo-Coullet-Tresser system}

In the paper \cite{ArnCoulTres1982}, a smooth three-dimensional system
\begin{equation}
\begin{cases}
\dot x = y \\
\dot y = z \\
\dot z = -y -\beta z + \mu x (1-x),
\end{cases}
\label{eq:ACTsys}
\end{equation}
was proposed which demonstrates spiral chaos at certain regions of values of the parameters $\beta$ and $\mu$. Let us consider this system in more detail.

Since system \eqref{eq:ACTsys} has the constant divergence equal to $-\beta$, attractors can exist only for $\beta>0$. Note that for all values of parameters $\mu> 0$ and $\beta> 0$, system \eqref{eq:ACTsys} has two equilibria  $O_1(0,0,0)$ and $O_2(1,0,0)$. Moreover, the equilibrium $O_1$ for $\beta <\sqrt{3}$ is always a saddle-focus (2,1), while the equilibrium $O_2$ can be both stable and saddle-focus (1,2).

\begin{figure}[!h]
\begin{minipage}[h]{0.32\linewidth}
\center{\includegraphics[width=1\linewidth]{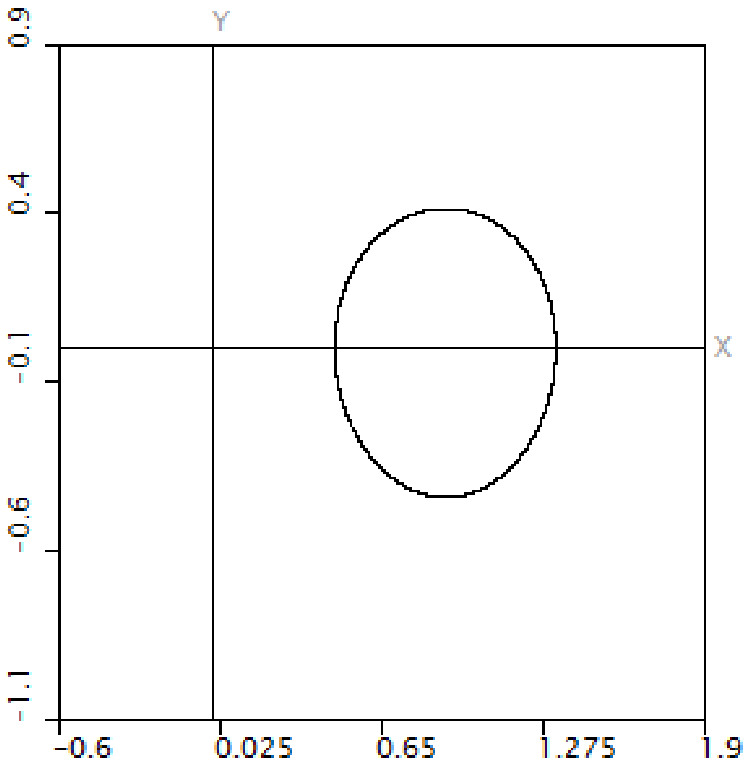} \\ a) {$\mu = 0.5$}}
\end{minipage}
\hfill
\begin{minipage}[h]{0.32\linewidth}
\center{\includegraphics[width=1\linewidth]{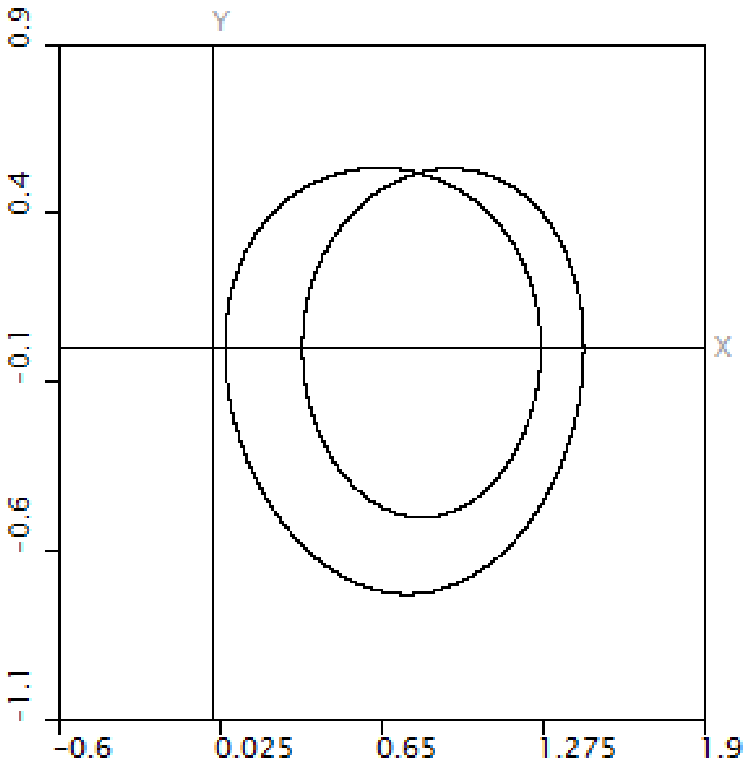} \\ b) {$\mu = 0.75$}}
\end{minipage}
\hfill
\begin{minipage}[h]{0.32\linewidth}
\center{\includegraphics[width=1\linewidth]{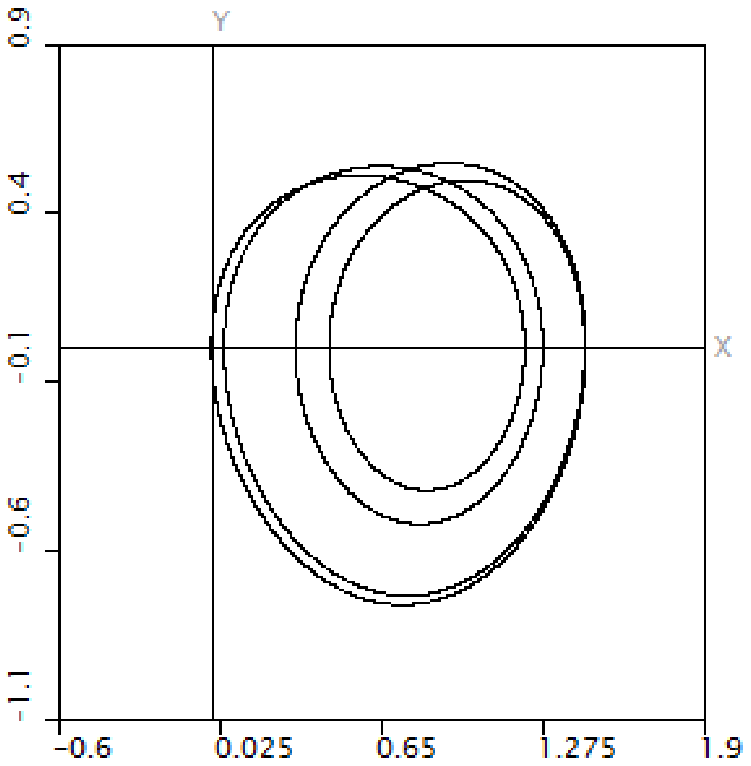} \\ c) {$\mu = 0.776$}}
\end{minipage}
\vfill
\begin{minipage}[h]{0.32\linewidth}
\center{\includegraphics[width=1\linewidth]{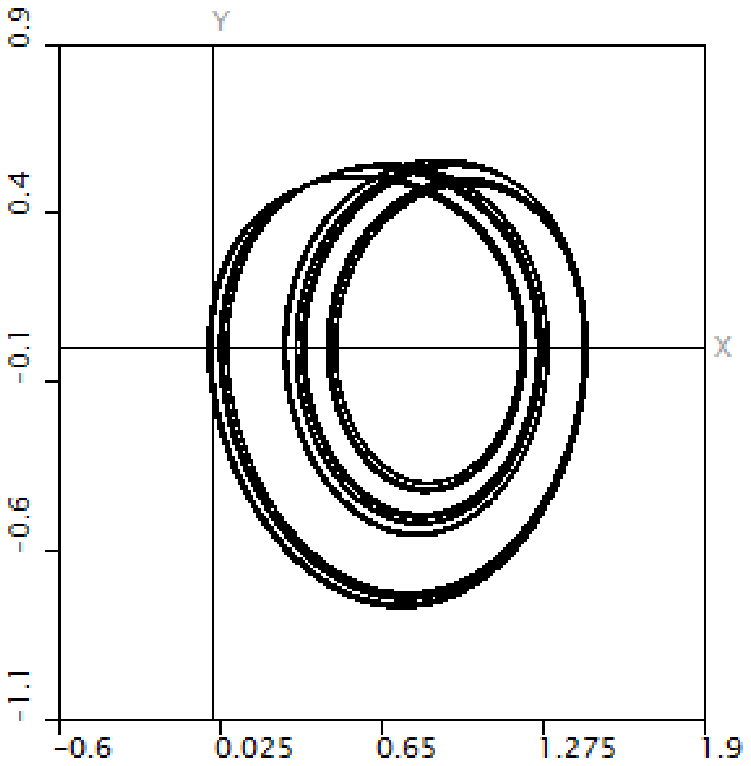} \\ d) {$\mu = 0.782$}}
\end{minipage}
\hfill
\begin{minipage}[h]{0.32\linewidth}
\center{\includegraphics[width=1\linewidth]{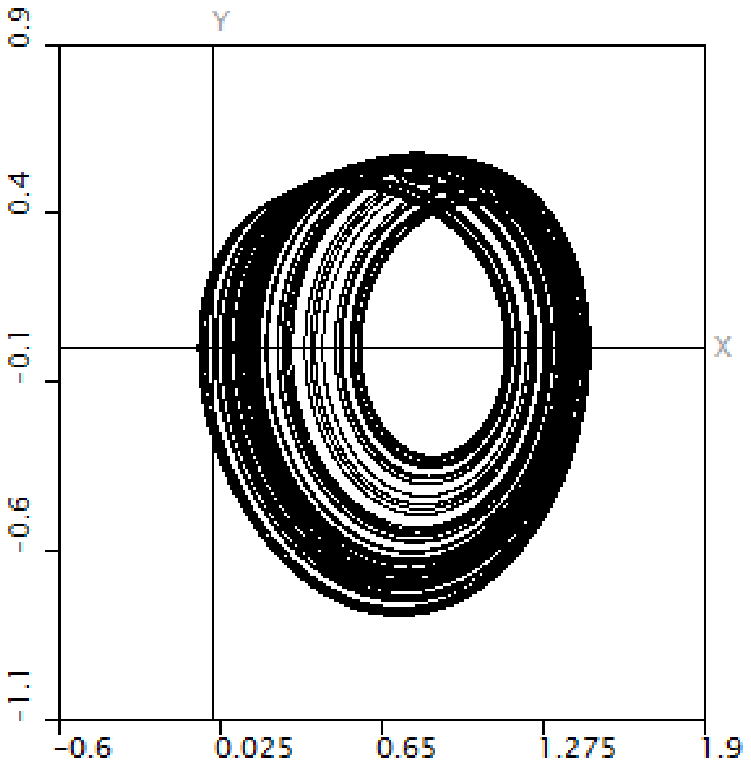} \\ e) {$\mu = 0.8$}}
\end{minipage}
\hfill
\begin{minipage}[h]{0.32\linewidth}
\center{\includegraphics[width=1\linewidth]{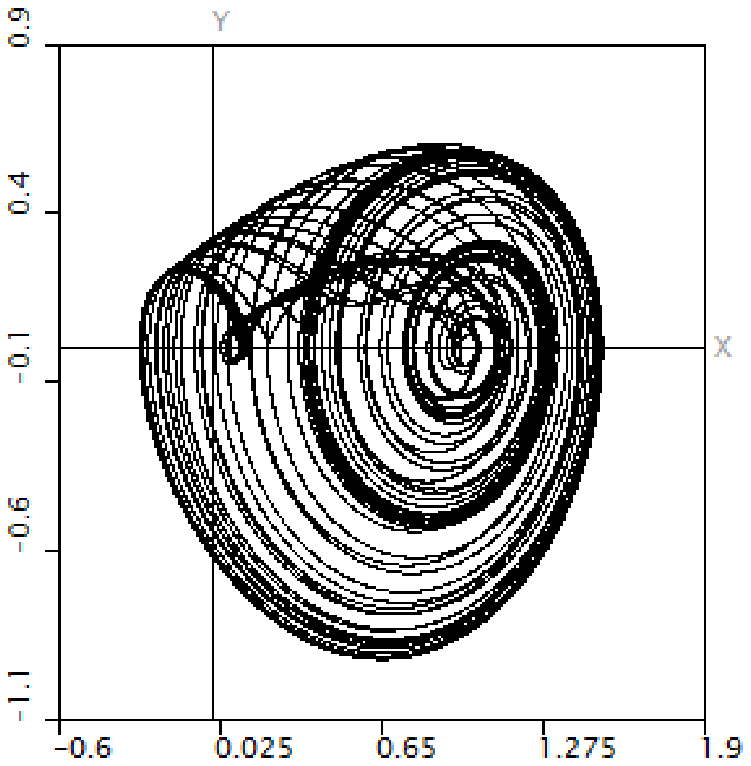} \\ f) {$\mu = 0.863$}}
\end{minipage}
\caption{{\footnotesize Scenario of the emergence of Shilnikov attractor in system \eqref{eq:ACTsys}.
}}
\label{fig:Arneodo_scenario}
\end{figure}

Following \cite{ArnCoulTres1982} we put $\beta = 0.4$. Now our goal is to illustrate the scenario of the emergence of the Shilnikov attractor containing the saddle-focus $O_2$ when the parameter $\mu$ changes. In Fig.~\ref{fig:Arneodo_scenario} main stages of the scenario are shown.

For $0<\mu < \mu_1 = 0.4$, the stable equilibrium $O_2$ is the only attractor of the system. Note that the boundary of its absorbing domain is formed by the two-dimensional stable invariant manifold of the saddle-focus $O_1$. At $\mu = \mu_1$, the equilibrium $O_2$ undergoes a supercritical Andronov-Hopf bifurcation. As a result, for $\mu > \mu_1$, a stable limit cycle $L$ is born in a neighborhood of the equilibrium $O_2$ that becomes a saddle-focus (1,2), see Fig.~\ref{fig:Arneodo_scenario}a. The stable limit cycle exists for $\mu_1 < \mu < \mu_2 \approx 0.72$. Starting with $\mu = \mu_2$, the cycle goes through a cascade of period-doubling bifurcations (see Fig.~\ref{fig:Arneodo_scenario}b,c after the first and second period-doubling bifurcations), and, as a result, a strange attractor of Feigenbaum type
appears, see Fig.~\ref{fig:Arneodo_scenario}d. Then, this attractor is transformed into a R\"ossler-like attractor (on a two-dimensional Poincar\'e section, it can be also interpreted as a H\'enon-like attractor).
With a further increase in the parameter $\mu$ up to $\mu = \mu_3 \approx 0.86311445$, the attractor changes quite smoothly
and is not homoclinic, since it is separated from the saddle-focus $O_2$.

\begin{figure}[!ht]
\begin{minipage}[h]{0.32\linewidth}
\center{\includegraphics[width=1\linewidth]{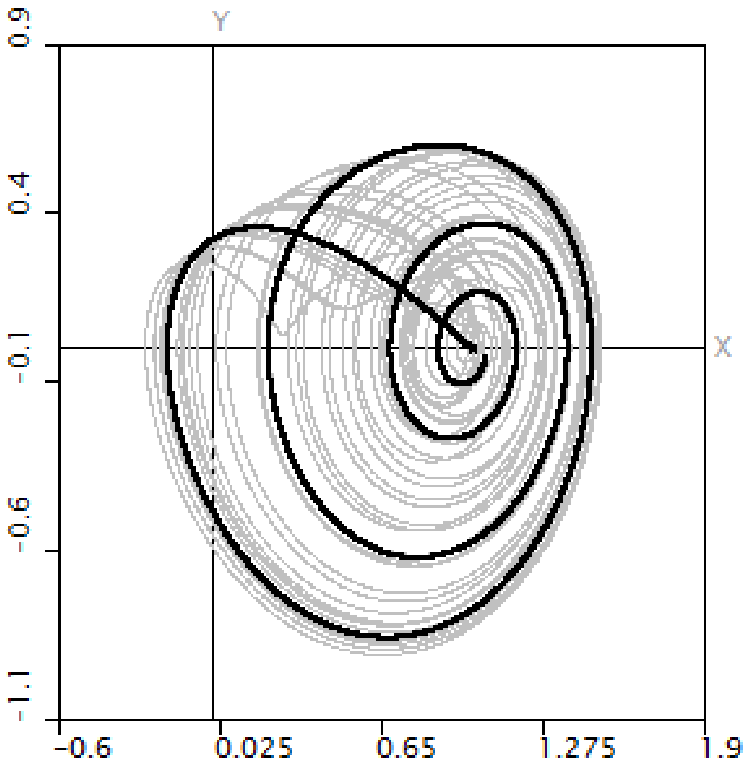} \\ a) {$\mu = 0.86311445$}}
\end{minipage}
\hfill
\begin{minipage}[h]{0.32\linewidth}
\center{\includegraphics[width=1\linewidth]{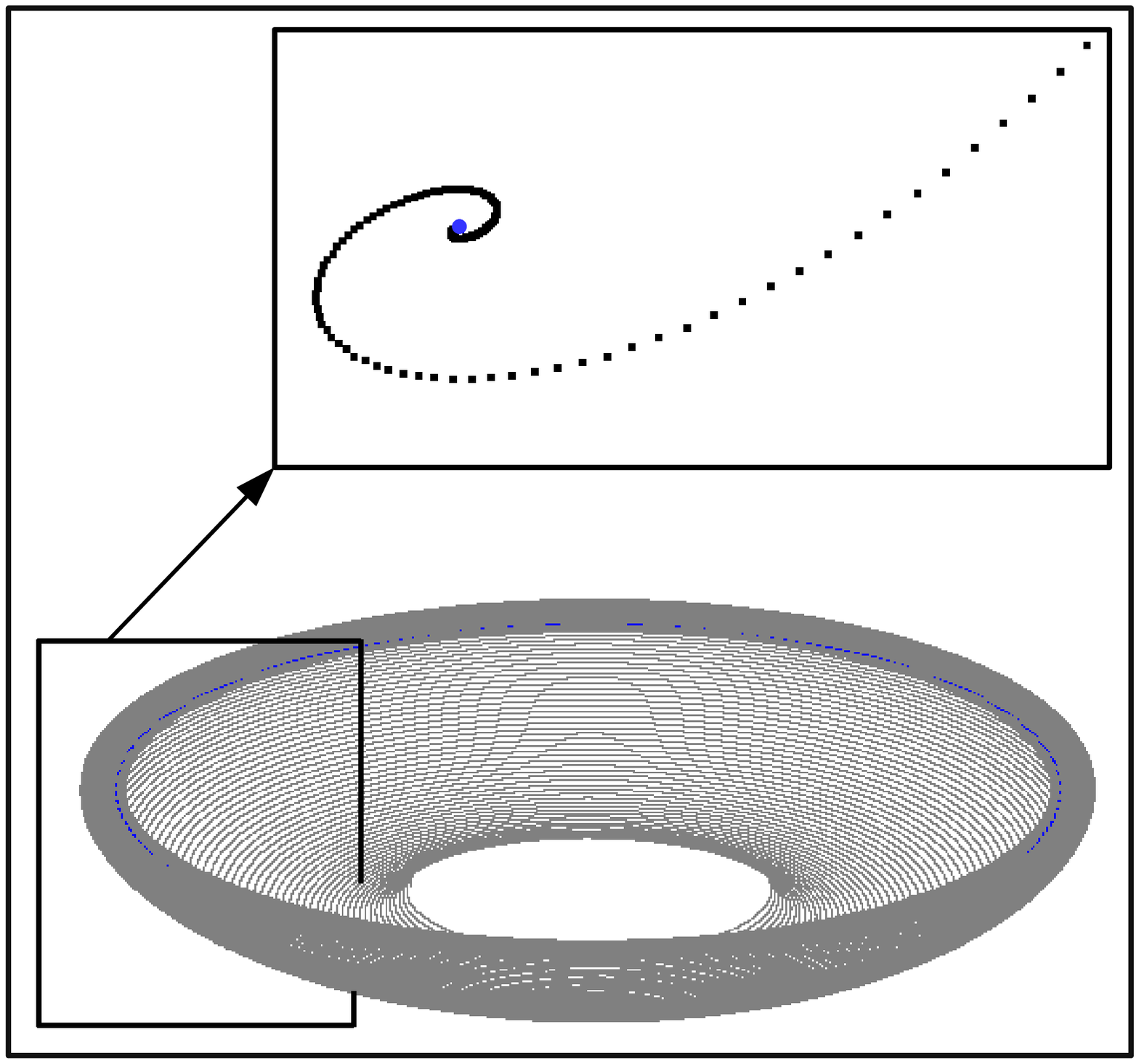} \\ b) {$\beta = 0.01$, $\mu = 0.02$}}
\end{minipage}
\hfill
\begin{minipage}[h]{0.32\linewidth}
\center{\includegraphics[width=1\linewidth]{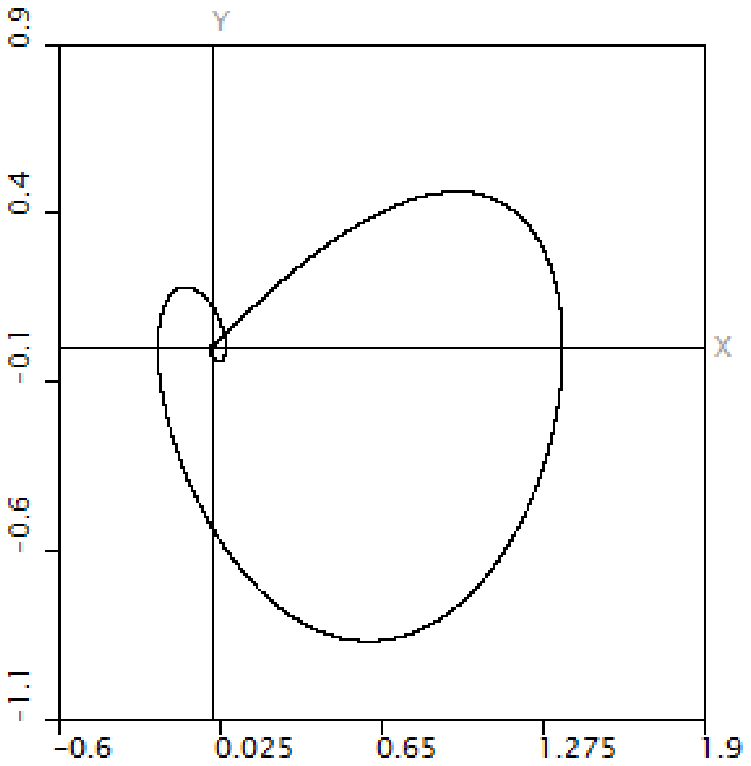} \\ c) {$\mu = 1.6062$}}
\end{minipage}
\caption{{\footnotesize (a) Shilnikov attractor with a homoclinic loop of the saddle-focus $O_2$ of type (1,2) (b) illustration of the formation of the Shilnikov whirlpool; (c) $\mu = \mu_5$, a homoclinic loop of the saddle-focus $O_1$ of type (2,1), here the attractor no longer exists.}}
\label{fig:Arneodo_loop2}
\end{figure}

However, as $\mu \to \mu_3$, a  ``hole'' around $O_2$ decreases and, finally, disappears at $\mu=\mu_3$.  In this case, orbits of the attractor can come arbitrarily close to the saddle-focus $O_2$, see Fig.~\ref{fig:Arneodo_scenario}f.
As we know, this is due to the appearance of a homoclinic orbit to the saddle-focus $O_2$ and, hence, to the creation of
the Shilnikov attractor, see  Fig.~\ref{fig:Arneodo_loop2}a, where the numerically found homoclinic orbit is also shown.

Note that, using system \eqref{eq:ACTsys}, we can trace a formation of the Shilnikov whirlpool -- when the cycle $L$ becomes of focal type. Illustration of the moment of occurrence of the Shilnikov whirlpool is shown in Fig.~\ref{fig:Arneodo_loop2}b, where, for greater clarity, we take  $\beta = 0.01$ (then  the divergence of the system $Div = -\beta$ will be small) in order to be able to show how the unstable two-dimensional manifold of the saddle-focus $O_2$ winds around the stable limit cycle $L$.

With a further increase in the parameter $\mu$, from $\mu=\mu_3$ to $\mu = \mu_{4} \approx 0.873$, a phenomenon of complication of the structure of the Shilnikov whirlpool can be observed, see e.g. \cite{LetDutMah1995, Malykh20}. At $\mu = \mu_{4}$, a crisis of the attractor occurs:  it collides with the two-dimensional stable manifold of the saddle-focus $O_1$ which is a natural boundary of the absorbing domain for the attractor under consideration. Note that the saddle-focus $O_1$ can also have homoclinic loops, see Fig.~\ref{fig:Arneodo_loop2}a where such a loop is shown for   $\mu = \mu_5 \approx 1.6062$ (it was also found numerically in \cite{ArnCoulTres1982}). However, any attractors do not exist for close values of $\mu$.

\subsection{Example of Shilnikov attractor in the Gaspard-Nicolis model of chemical oscillator.}

In this section  we consider an example of interesting system in which the Shilnikov attractor arises almost instantly as a result of a single rigid (subcritical) bifurcation. This scenario begins as usual:  a stable equilibrium $O$ undergoes a supercritical Andronov-Hopf bifurcation, as a result, it becomes a saddle-focus (1,2) and a stable limit cycle $L$ is born. This cycle first becomes focal one and the Shilnikov whirlpool (with the boundary $W^u(O)$) is formed. However, this limit cycle itself remains stable for a long time (by the parameter) and, independently of it,  non-attractive chaotic dynamics (metastable chaos)  has time to develop  inside the whirlpool. Next, at a certain value of the parameter, the cycle $L$ undergoes a period doubling  bifurcation that is subcritical:  a saddle cycle of doubled period merges with $L$. After this, a strange attractor is immediately observed, which at very close values of a parameter becomes homoclinic -- the Shilnikov attractor appears.

Below we will show how this happens in the system
\begin{equation}
\begin{cases}
\dot x = x(\beta x - f y - z + g), \\
\dot y = y (x + s z - \alpha), \\
\dot z = (x - \alpha z^3 + b z^2 - c z) / \varepsilon ,
\end{cases}
\label{eq:COsys}
\end{equation}
proposed by P. Gaspard and G. Nicolis in the paper \cite{GaspardNicolis83}, in which it was shown that  strange attractors containing a saddle-focus (1,2) can exist in system (\ref{eq:COsys}). The corresponding parameter regions where such attractors can be observed have been specified in \cite{Gallas2010}. Following \cite{Gallas2010} we take some parameters to be fixed

\begin{equation}
b = 3, \varepsilon = 0.01, f = 0.5, g = 0.6, s = 0.3, c = 4.8, \alpha = 0.7825,
\label{eq:COparams}
\end{equation}
and the parameter $\beta$ choose as the governing one.

\begin{figure}[!ht]
\begin{minipage}[h]{0.32\linewidth}
\center{\includegraphics[width=1\linewidth]{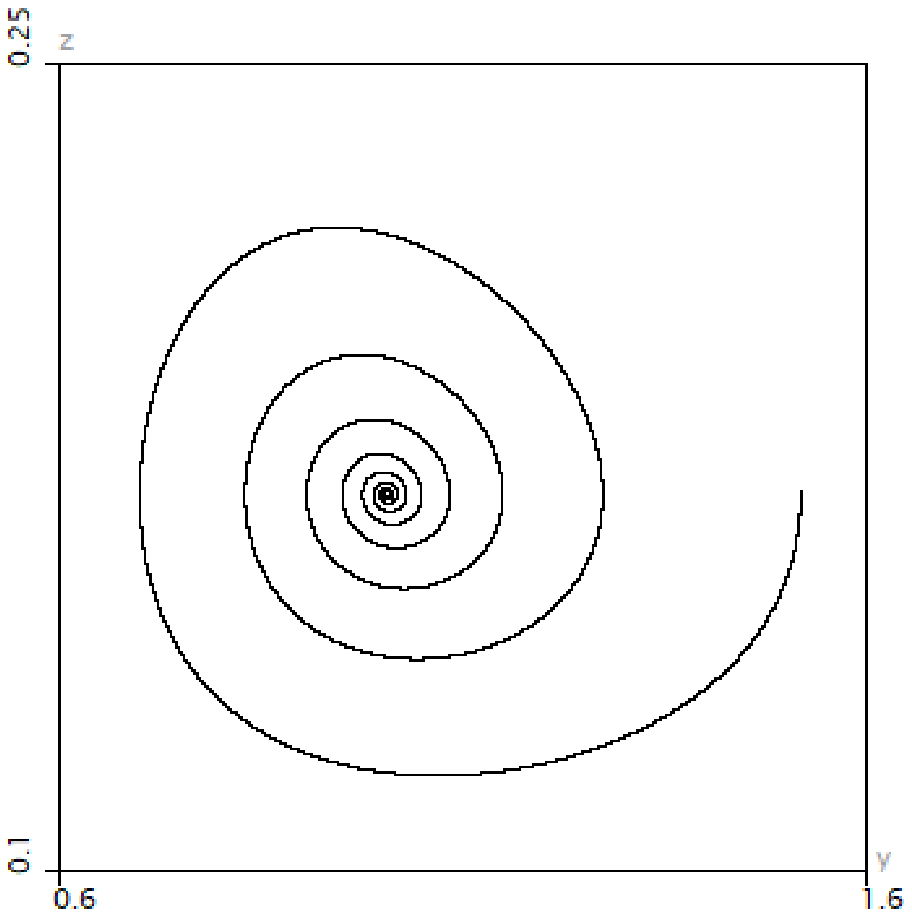} \\ a) {$\beta = 0.1$}}
\end{minipage}
\hfill
\begin{minipage}[h]{0.32\linewidth}
\center{\includegraphics[width=1\linewidth]{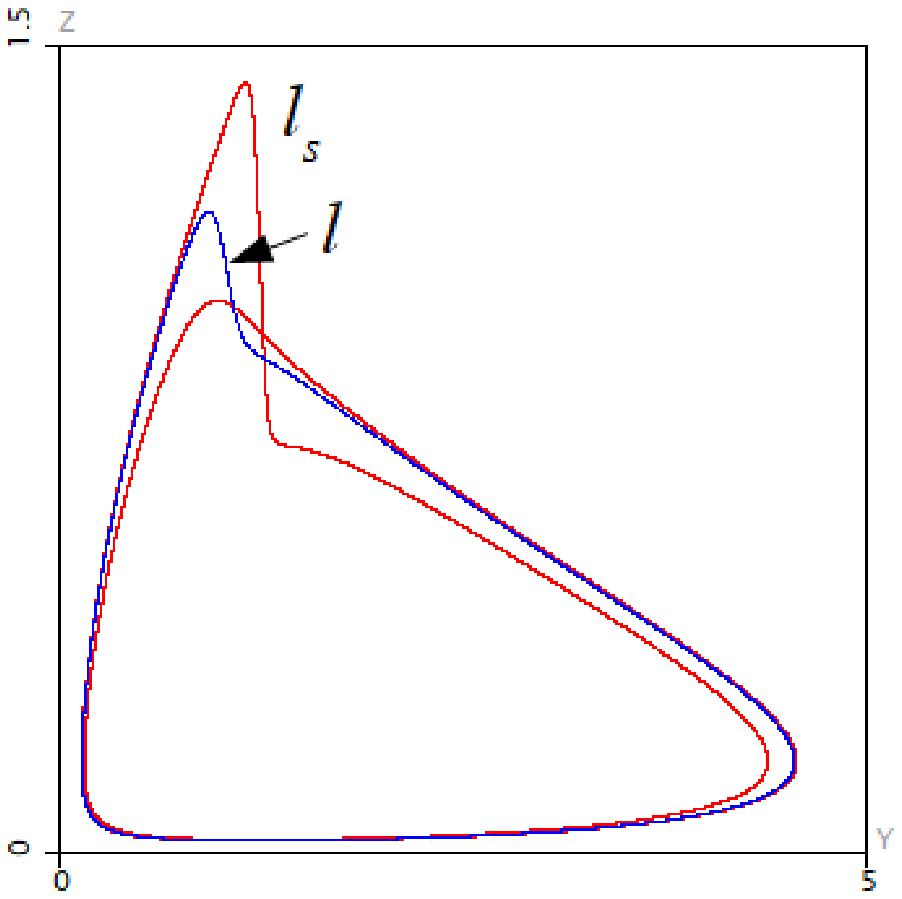} \\ b) {$\beta = 0.3812$}}
\end{minipage}
\hfill
\begin{minipage}[h]{0.32\linewidth}
\center{\includegraphics[width=1\linewidth]{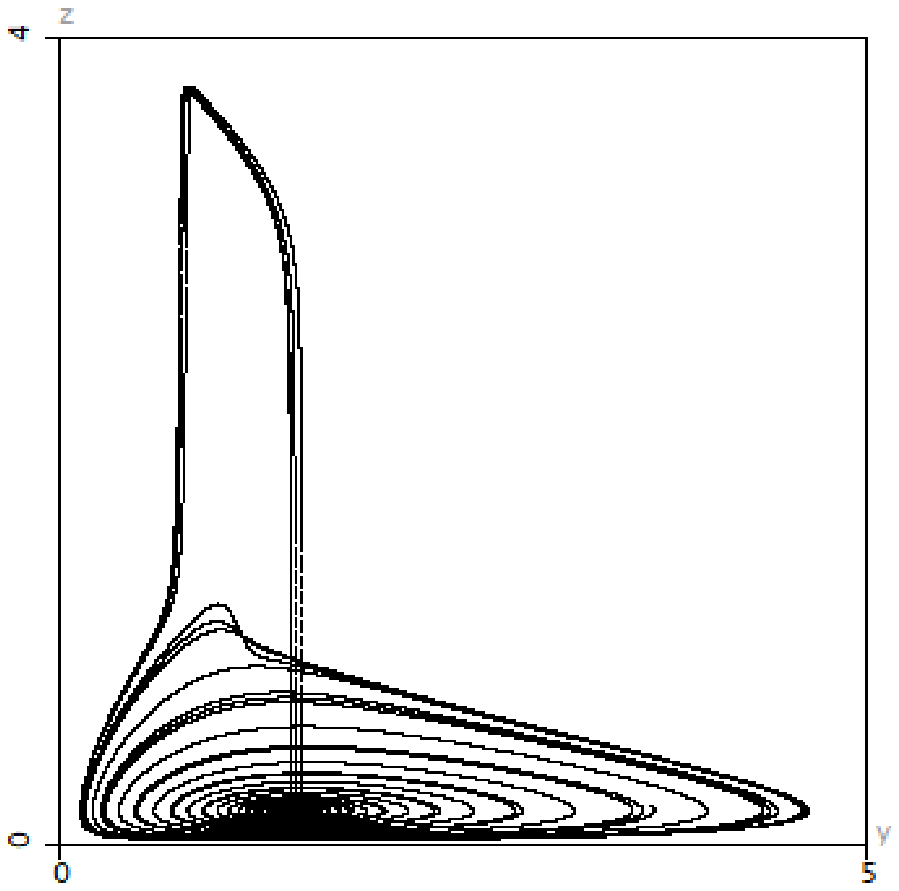} \\ c) {$\beta = 0.39213$}}
\end{minipage}
\caption{{\footnotesize Scenario of the Shilnikov attractor emergence in the model \eqref{eq:COsys} of the chemical oscillator: (a) the attractor is the stable equilibrium state $O$; (b) saddle cycle $l_s$ of period-2 approaches the stable limit cycle $l$; (c) Shilnikov homoclinic attractor.}}
\label{fig:CO_scenario}
\end{figure}

For $\beta < \beta_{AH} \approx 0.261$, the nonzero equilibrium $O$ of system  (\ref{eq:COsys}) is stable, see Fig.~\ref{fig:CO_scenario}a ($O$ is only equilibrium with positive coordinates having a physical sense). For $\beta = \beta_{AH}$ this equilibrium undergoes a supercritical Andronov-Hopf bifurcation. After this, on the interval $\beta \in (\beta_{AH}, \beta_{PD} \approx 0.3817)$, the attractor of the system is the stable limit cycle $l$, see Fig.~\ref{fig:CO_scenario}b, and the equilibrium $O$ becomes a saddle-focus (1,2).
At $\beta \approx 0.365$ the cycle $l$ becomes focal and  a Shilnikov funnel is formed, further  both multipliers of $l$ becomes negative and one of them tends to $-1$ as
$\beta \rightarrow \beta_{PD}$. Simultaneously, a saddle limit cycle
$l_s$ of double period approaches the cycle $l$, see Fig.~\ref{fig:CO_scenario}b and merges with $l$. As a result of this {\em subcritical} period doubling bifurcation, at $\beta> \beta_{PD}$ the limit cycle $l$ becomes saddle  and, instantly, a strange attractor is observed.  With a further increase in the parameter $\beta$, the orbits of this attractor begin to approach closer and closer to the saddle-focus $O$. At $\beta = \beta_h \approx 0.3921$, a homoclinic loop of the saddle-focus $O$ is formed, i.e. a homoclinic Shilnikov attractor arises, Fig.~\ref{fig:CO_scenario} c.

\section{On Shilnikov scenarios for three-dimensional orientable and nonorientable maps} \label{sec:Shiscen_map}

We discuss now scenarios of the emergence of  \emph{discrete Shilnikov attractors}, i.e. homoclinic attractors containing a fixed (periodic) point that is a saddle-focus with two-dimensional unstable manifold, thus, it is a saddle-focus (1,2) in the case of three-dimensional maps.

We consider one-parameter families $T_\mu$ of three-dimensional diffeomorphisms in two cases, when map $T_\mu$ is orientable (orientation preserving) and when $T_\mu$ is nonorientable (orientation reversing) map. In order not to get involved with the problems of orientability of the ambient  manifold, we will assume that in both cases $T_\mu$ is given in $R^3$. Then the Jacobian $J(T_\mu)$ of the map $T_\mu$ will be everywhere positive in the orientable case and everywhere negative in the nonorientable case.

\subsection{The orientable case} \label{sec:Shor_scen}

\begin{figure}[ht]
\centerline{\epsfig{file=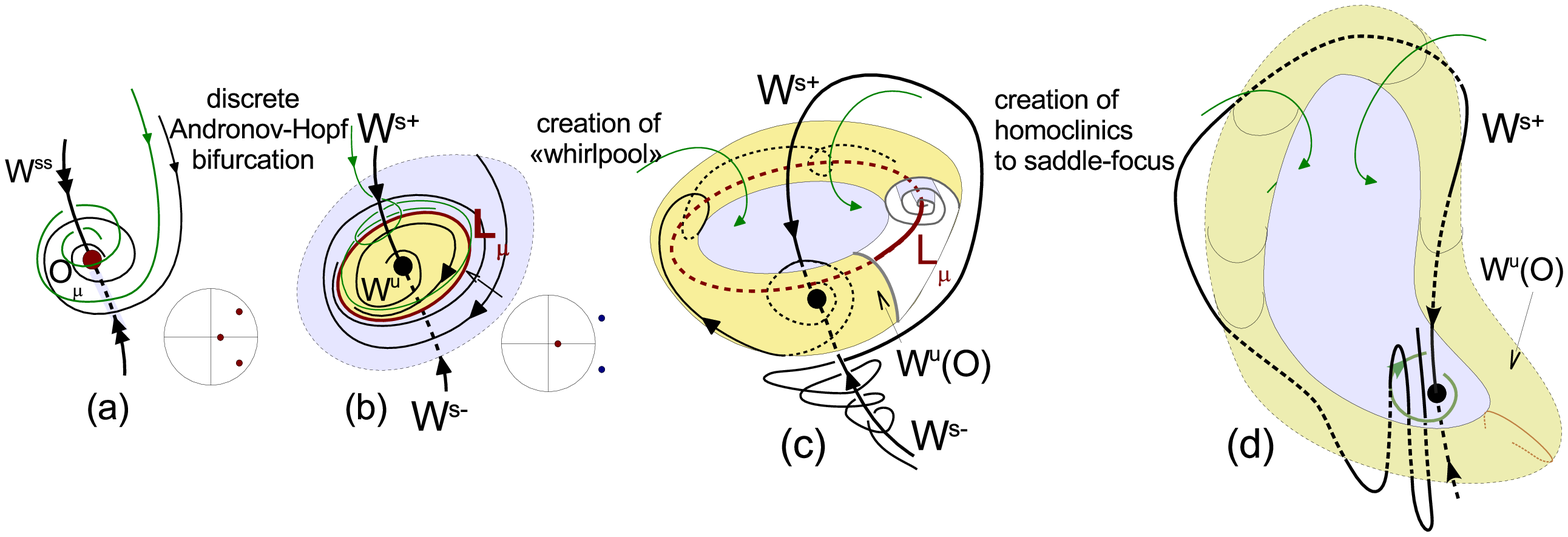, width=14cm
}}
\vspace{-1cm}
\caption{{\footnotesize A sketch of scenario of the emergence of a discrete Shilnikov attractor.}}
\label{shildiscr}
\end{figure}

A sketch of scenario of a typical discrete Shilnikov attractor appearance for one-parameter families $T_\mu$ of three-dimensional orientable maps is illustrated in Fig.~\ref{shildiscr}.
This scenario starts with a stable fixed point $O_\mu$ that loses the stability at $\mu=\mu_1$ under a discrete supercritical (soft) Andronov-Hopf bifurcation,
Fig.~\ref{shildiscr}a--b:  for
$\mu>\mu_1$ the point $O_\mu$ becomes a saddle-focus (1,2)
and a stable closed invariant curve $L_\mu$ is born in a neighborhood of $O_\mu$. Thus,
the curve $L_\mu$ becomes the attractor of map $T_\mu$.
During this transition, see Fig.~\ref{shildiscr}a--b, the point $O_\mu$ changes its type from a stable point to a saddle-focus (1,2) point:  for $\mu <\mu_1$ it has three multipliers inside the unit circle, for $\mu=\mu_1$ two complex conjugate multipliers of $O_\mu$ fall on the unit circle, and for $\mu >\mu_1$ they go outside it.

The next stage of the scenario is connected with changes in $L_\mu$. Typically, this happens in the following way.  Just after the discrete Andronov-Hopf bifurcation, at small
$\mu-\mu_1$,  the unstable manifold of $O_\mu$ is a two-dimensional disk $D_\mu$ with the boundary $L_\mu$ that has a type of nodal invariant curve. We assume that, at further changing $\mu$, the curve $L_\mu$ undergoes a ``smooth bifurcation'', when it changes its type from nodal to focal and, thus, the two-dimensional manifold $W^u(O_\mu)$ begins to wind up on $L_\mu$, like a roll. As a result, $W^u(O_\mu)$ takes the form of a kind of cauldron, the Shilnikov whirlpool,  inside which all orbits from the absorbing region are drawn, except for orbits of the stable separatrix $W^{s-}(O_\mu)$,
see Fig.~\ref{shildiscr}c.

Then, chaotic dynamics begins to develop in this whirlpool as $\mu$ changes. At first, the attractor is simple,  it is the stable invariant curve $L_\mu$, and then it loses its stability. This can happen in a variety of ways (for example, in a soft way through a cascade of doubling of invariant curves with their subsequent destruction and the formation of attractors of torus-chaos type, or in a rigid way through a subcritical bifurcation with some of stable invariant curves, after which chaos can be observed immediately,
``by explosion'', etc.).
In principle, for the essence of the phenomenological scenario, it does not matter how this happens. The main thing here is that the invariant manifolds
$W^u(O_\mu)$ and $W^{s+}(O_\mu)$  begin to intersect and a strange homoclinic attractor can arise containing the saddle-focus
$O_\mu$ and entirely its two-dimensional unstable manifold, see Fig.~\ref{shildiscr}d.
We call this attractor a {\em discrete Shilnikov attractor}.

It should be noted that there is almost complete similarity in the main features of this scenario with the corresponding scenario in the case of a flow, see section 1. However, even here one can see a significant difference. The Shilnikov homoclinic flow attractor exists only for discrete values of the control parameter corresponding to the existence of homoclinic loops of the saddle-focus equilibrium. Whereas, the Shilnikov attractor for maps exists on  intervals of parameter values for which the invariant manifolds ($W^u(O_\mu)$ and $W^{s+}(O_\mu)$) have transversal intersections. Moreover, these intervals can be large enough and their values can reach even those at which the attractor is destroyed and disappears altogether.

\begin{figure}[!ht]
\begin{minipage}[h]{0.49\linewidth}
\center{\includegraphics[width=1\linewidth]{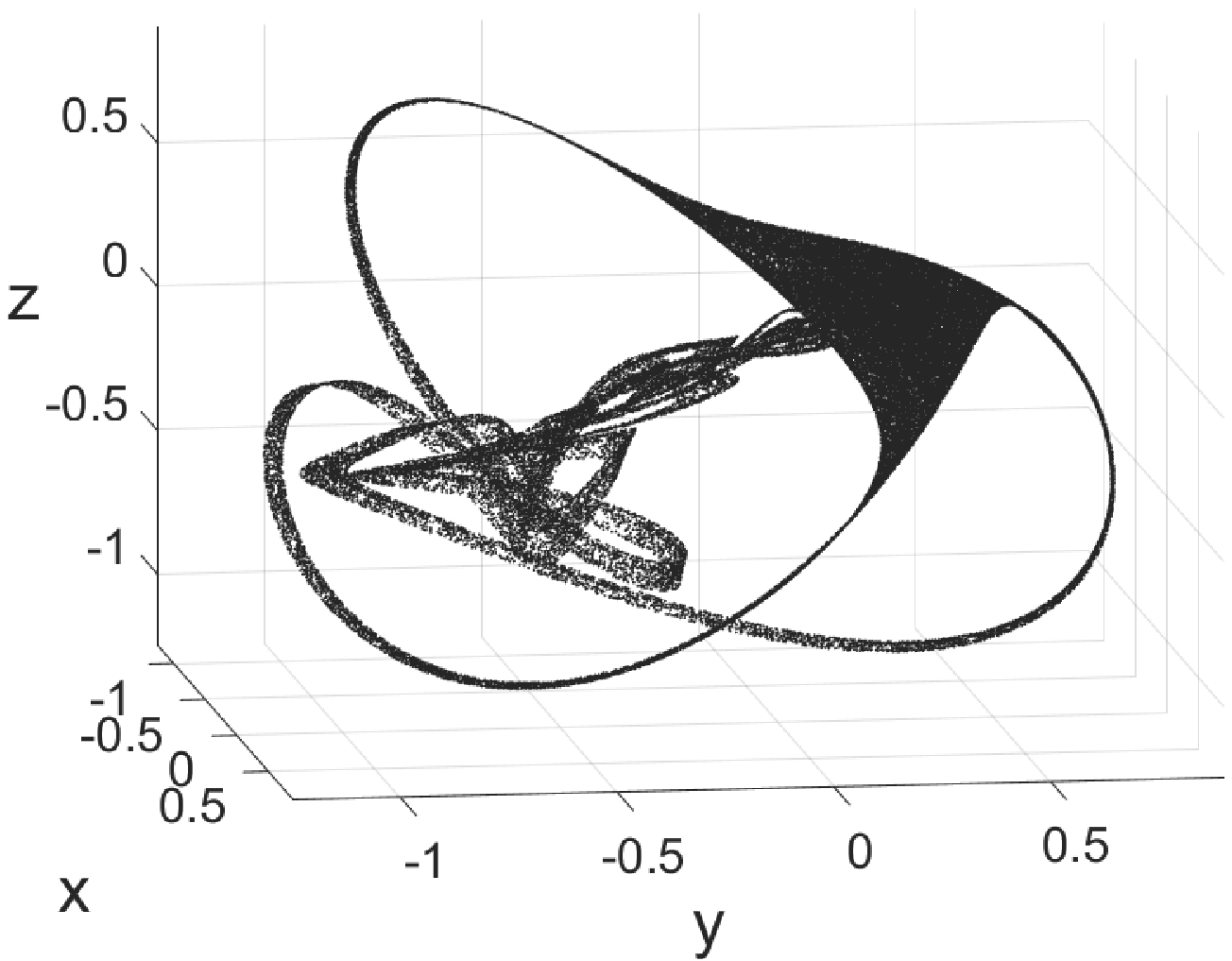} \\ a)}
\end{minipage}
\hfill
\begin{minipage}[h]{0.49\linewidth}
\center{\includegraphics[width=1\linewidth]{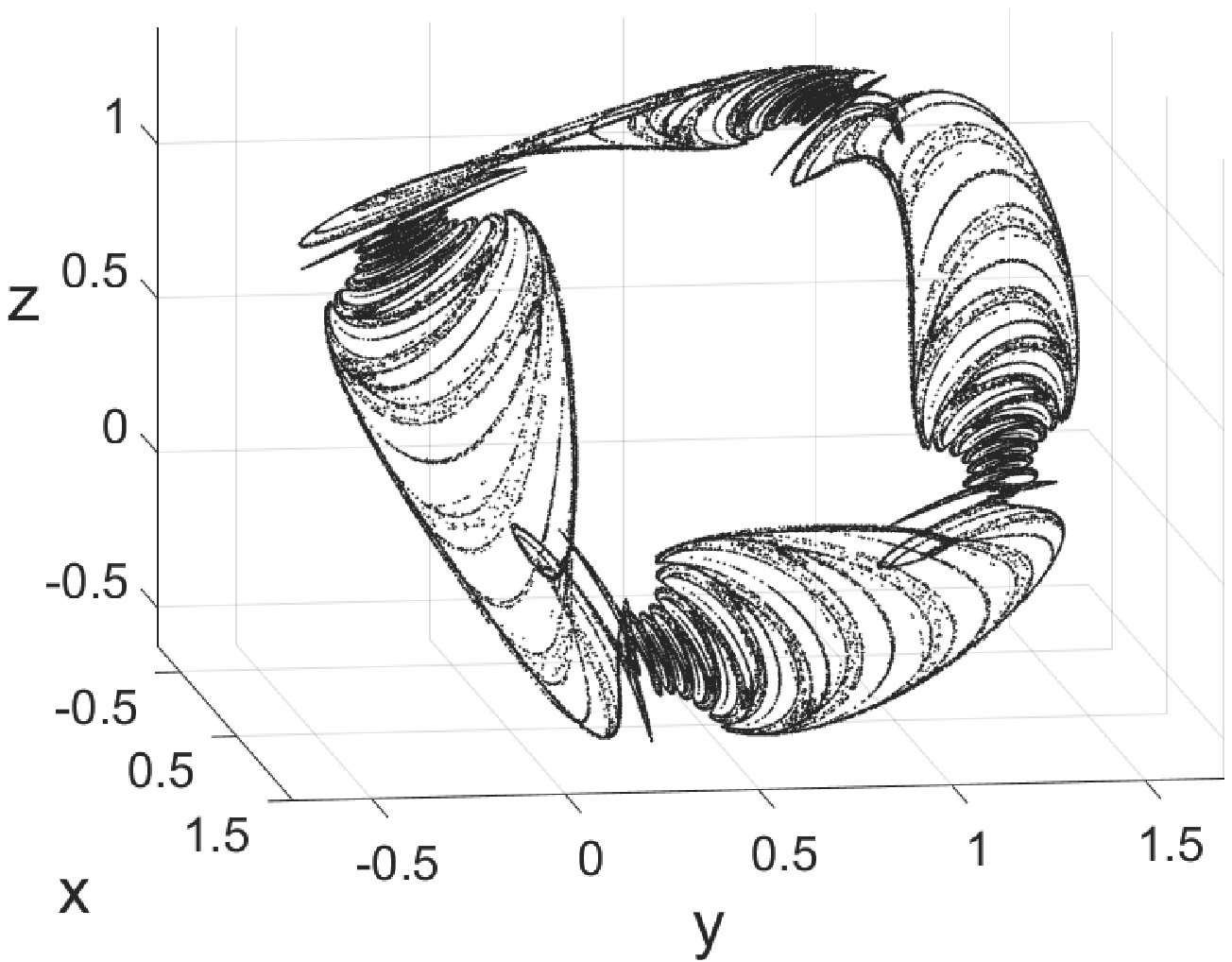} \\ b)}
\end{minipage}
\caption{{\footnotesize Discrete attractors in resonant case: (a) a ``triangular'' discrete Shilnikov attractor; (b) a ``superspiral'' period-4 attractor. We found these attractors in the three-dimensional Mir\'a map
$\bar x = y, \bar y = z, \bar z = M_1 + B x + M_2 z - y^2$ for $B=0.7$ and (a) $M_1 = 0.195, M_2 = -0.26$; (b) $M_1 = 0.35, M_2 = 0.8$.}}
\label{Shattr_rez}
\end{figure}

There is also one more important feature of Shilnikov discrete attractors, which flow attractors do not have. This feature manifests itself in the case when the stable invariant curve $L_\mu$ is resonant. Then, on the curve itself, there are alternating saddle and stable periodic points of the same period, and the formation of the Shilnikov whirlpool occurs due to the fact that the stable points become foci. Thus, the manifold $W^u(O_\mu)$ is twisted over $L_\mu $ only in some places (near stable points). This also subsequently affects the shape of the emerging Shilnikov attractor. For example, such an attractor can have a characteristic ``triangular '' or ``square'' shape in the case of resonances 1:3 and
1:4, when $O_\mu$ has a pair of multipliers $\lambda e^{\pm i\psi}$ with $\psi$ close to $2\pi/3$ and $\pi/2$,  respectively, see Fig.~\ref{Shattr_rez}. In addition, the resonant invariant curves of three-dimensional maps can themselves be destroyed in very interesting ways, giving rise to amusing attractors inside the funnel, see, for example, Fig.~\ref{Shattr_rez}b, which shows a ``superspiral'' attractor containing two orbits of period 4, which are saddle-foci of type (2,1) and (1,2).

\subsection{The nonorientable case}  \label{sec:Shnor_scen}

\begin{figure}[ht]
\centerline{\epsfig{file=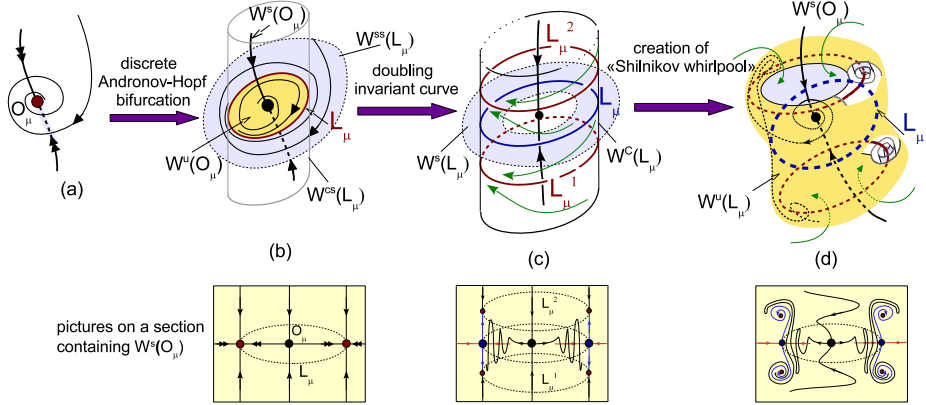, width=14cm
}}
\vspace{-1cm}
\caption{{\footnotesize A sketch of scenario of the emergence of a nonorientable discrete Shilnikov attractor.}}
\label{scen-shil-nor4}
\end{figure}

A sketch of a typical discrete Shilnikov scenario for one-parameter families $T_\mu$ of three-dimensional nonorientable maps is illustrated in Fig.~\ref{scen-shil-nor4}. This scenario starts with those $\mu$ at which $T_\mu$ has a nonorientable stable fixed point $O_\mu$ with multipliers $-\lambda, \gamma^{\pm i\psi}$, where $0<\lambda<\gamma<1$ and $0<\psi<\pi$, Fig.~\ref{scen-shil-nor4}a. We assume that $O_\mu$ loses the stability at $\mu=\mu_1$ under a supercritical Andronov-Hopf bifurcation: for $\mu>\mu_1$ the point $O_\mu$ becomes a nonorientable saddle-focus (1,2), and a stable closed invariant curve $L_\mu$ is born in a neighborhood of $O_\mu$, Fig.~\ref{scen-shil-nor4}b. Thus, the curve $L_\mu$ becomes the attractor of map $T_\mu$. During this transition, the point $O_\mu$ has three multipliers less than one in the absolute value for $\mu <\mu_1$, then, for $\mu=\mu_1$, two of its complex conjugate multipliers fall on the unit circle, and for $\mu >\mu_1$ they go out and the point $O_\mu$ becomes a saddle-focus (1,2) that is nonorientable since its stable multiplier is negative.

Next stage of the scenario is connected with changes in $L_\mu$. Typically, this happens in the following way.  Just after the Andronov-Hopf bifurcation, at small $\mu-\mu_1$,  the unstable manifold of $O_\mu$ is a two-dimensional disk $D_\mu$ with the boundary $L_\mu$ where the curve $L_{\mu}$ has a type of nonorientable nodal invariant curve. Because of the curve $L_\mu$ is nonorientable it can not become of focal type as in the orientable case. Thus, the manifold $W^u(O_\mu)$ can not take a form of whirlpool just over $L_\mu$. However, there exists another way for creation of a whirlpool. Namely, at further changing $\mu$, first the curve $L_\mu$ undergoes a doubling bifurcation: the curve $L_\mu$ itself becomes saddle and two stable period-2 invariant curves $L_\mu^1$ and $L_\mu^2$  originate from it here $L_\mu^2 = T_\mu(L_\mu^1)$ and $L_\mu^1 = T_\mu(L_\mu^2)$. Such a feature of the bifurcation of doubling of an invariant curve is obtained due to the fact that near the bifurcation the curve $L_\mu$ possesses two invariant manifolds,  strongly stable $W^{ss}$ and central $W^c$, which locally are both cylinders, but the map $T_\mu$ in the restriction to $W^c$ changes orientation. Accordingly, after the bifurcation, the curves $L_\mu^1$ and $L_\mu^2$ lie on $W^c$ on opposite sides of $L_\mu$.\footnote{More details about bifurcations of doubling of closed invariant curves can be found, for example, in \cite{GGT21}.} After this doubling bifurcation, the unstable manifold $W^u(O_\mu)$ immediately corrugates and begins to rush between the curves $L_\mu^1$ and $L_\mu^2$ but still remains inside the cylinder $W^c(L_\mu)$. Note that the curves $L_\mu^1$ and $L_\mu^2$ are both invariant and orientable for $T^2_\mu$, therefore they can become focal ones when changing $\mu$. In this case the manifold $W^u(O_\mu)$ starts to wind on both $L_\mu^1$ and $L_\mu^2$ and, thus, a nonorientable Shilnikov whirlpool is created, see Fig~\ref{scen-shil-nor4}d.

Then, in this whirlpool, when $\mu$ changes, chaotic dynamics begins to develop. At first, the attractor is simple, it is the stable period-2 invariant curve $(L_\mu^1,L_\mu^2)$, and then it loses its stability. Again, as for the orientable case, this can happen in a variety of ways... The main thing here is that the invariant manifolds $W^u(O_\mu)$ and $W^{s}(O_\mu)$  can intersect and a strange homoclinic attractor can arise containing the saddle-focus $O_\mu$ and its two-dimensional unstable manifold. In Fig~\ref{scen-shil-nor4}d below it is shown a sketch of manifolds $W^u(O_\mu)$ and $W^{s}(O_\mu)$ even before they crossed. One can imagine (but difficult to draw) what will happen when they intersect. However, it is clear that in this case the manifold $W^u(O_\mu)$, since the stable multiplier of the point $O_\mu$ is negative, will accumulate towards itself from both sides, and, accordingly, the point $O_\mu$ will reside inside the attractor. In the case of discrete orientable Shilnikov attractor, its fixed point $O_\mu$ lies on its boundary (since the global piece of $W^u(O_\mu)$ accumulates at $W^u_{loc}(O_\mu)$ only from one side (namely,  from the side where $W^{s+}_{loc}(O_\mu)$ is located).

Thus, we see that orientable and nonorientable discrete Shilnikov attractors have different structures. Although they both exist for an open set of parameter values, and it is also important for their geometry whether the invariant curve $L_\mu$ is resonant or not.

\section{Examples of three-dimensional maps with Shilnikov attractors} \label{sec:exShattr}

In this section we consider  two examples of discrete Shilnikov attractors in the
case of three-dimensional orientable and nonorientable maps. Recall that such attractors contain a
fixed point $O$ of the saddle-focus type, with eigenvalues $\lambda,
\gamma_{1,2} = \gamma e^{\pm i\psi}$, where $|\lambda|<1, \gamma >
1, 0<\psi<\pi$. The unstable manifold of $O$ is two-dimensional and it
resides entirely in the attractor.

\begin{figure}[ht]
    \center{\includegraphics[width=0.99\linewidth]{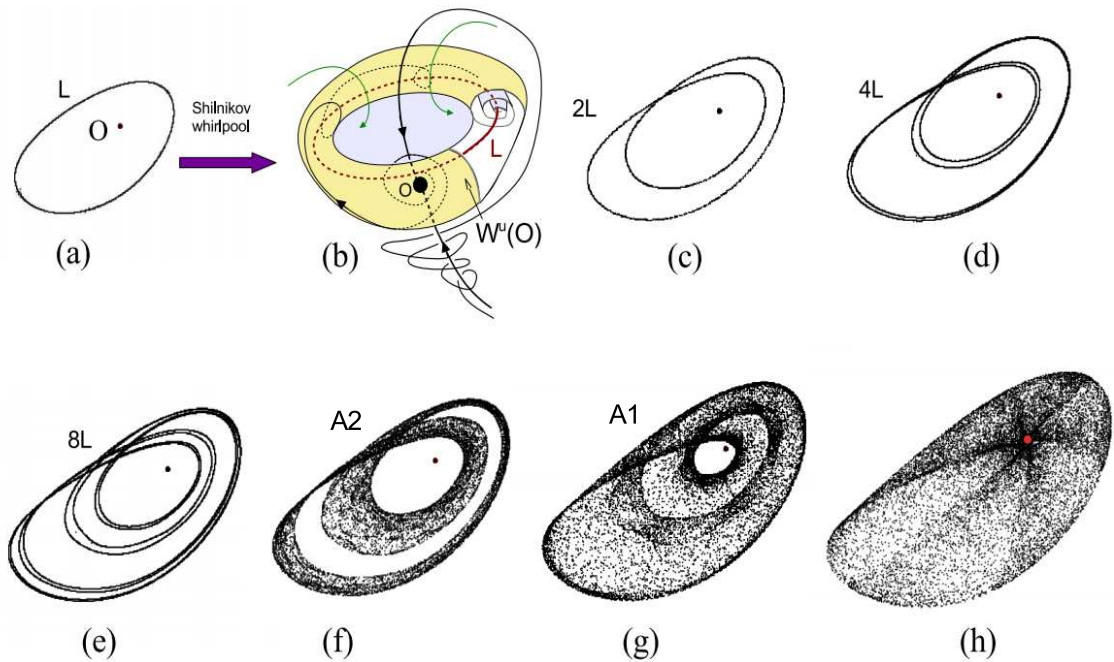}}
    \caption{{\footnotesize Stages of the emergence of the  Shilnikov attractor in
     map (\ref{Miramap}) for $B = 0.5, A = 1.49$, as $C$ changes (Figure (b) is only schematic):
     (a) $C = -1.7$, a stable invariant curve $L$ ; (b) a schematic picture of a Shilnikov funnel, which is formed when $W^u(O)$ begins
     to wound onto $L$ (the following figures (c)--(h) show what happens inside this funnel); (c) $C = -1.73$, the
     curve $L$ is doubled, $L\rightarrow 2L$; (d) $C = -1.76$ -- after the second doubling, $2L\rightarrow 4L$; (e) $C = -1.77$ --  after the third doubling, $4L\rightarrow 8L$; (f) $C = -1.775$,  a chaotic attractor $A2$ containing the saddle curves $4L$ and $2L$;
     (g) $C = -1.8$, a chaotic attractor $A1$ containing also the saddle curve $L$; (h) $C = -1.82$, a discrete Shilnikov attractor.}}
    \label{shil_or_mirmap}
\end{figure}

First, we consider the case of orientable
(with $0 <\lambda <1$)  Shilnikov attractor in three-dimensional map of the form
\begin{equation}
\bar x = y,\; \bar y = z, \; \bar z = Bx + Cy + Az-y^2,
\label{Miramap}
\end{equation}
where $A,B,C$ are parameters.
Numerically found example of such attractor is  shown in Fig.~\ref{shil_or_mirmap}h for map (\ref{Miramap}) with the Jacobian $B=0.5>0$.

Fig.~\ref{shil_or_mirmap} shows main stages of the development of a discrete Shilnikov  attractor in map (\ref{Miramap}) for $ B = 0.5,
A = 1.49 $ as $C$ changes. The formation of the attractor proceeds in accordance with the
\begin{figure}[ht]
\centerline{\centerline{\epsfig{file=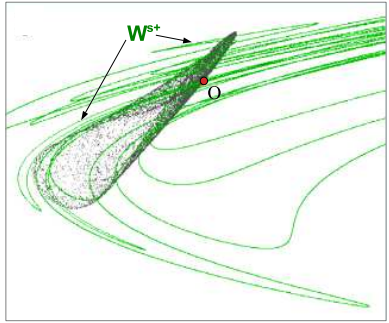, width=7cm}}}
\caption{{\footnotesize The Shilnikov attractor from Fig.~\ref{shil_or_mirmap}h and a part of the stable separatrix $W^{s+}(O)$ are shown in a suitable angle for viewing. }}
\label{shil_or_homoclinic}
\end{figure}
general bifurcation scenario described in
\cite{GGS12,GGKT14}, see section~\ref{sec:Shor_scen}:
\begin{itemize}
\item
At the beginning, the attractor of map \eqref{Miramap} is a stable fixed point $O(0,0,0)$ (it is stable
for $ -0.99> C> c_1 = - 1.495 $).
\item
Then, at $C = c_1$, the point $O$ undergoes a supercritical
Andronov-Hopf bifurcation: $O$ becomes a saddle-focus (1,2) and a stable closed invariant
curve $L$ is born in its neighborhood (the curve $L$ is shown in
Fig.~\ref{shil_or_mirmap}a at $ C= - 1.7$).
\item
Next, the unstable manifold of $O$ starts winding onto $L$ and the Shilnikov whirlpool is formed (schematically, the whirlpool is shown in
Fig.~\ref{shil_or_mirmap}b).
\item
Further, the dynamics inside the whirlpool become more and more
complicated.  In the case of  map (\ref{Miramap}), an example of
such development is shown in Fig.~\ref{shil_or_mirmap}d--h. First, we
observe three successive doubling bifurcations of stable invariant
curves, $L \to 2L \to 4L \to 8L$, Fig.~\ref{shil_or_mirmap}c-d-e.
We did not observe the doubling of $8L$. Instead,
a strange attractor appears which has at the beginning a torus-chaos type due to a destruction of the curve $8L$, then it sequentially captures the saddle curves $4L$, $2L$ (Fig.~\ref{shil_or_mirmap}f) and $L$ (Fig.~\ref{shil_or_mirmap}g).
\item
Finally, when homoclinic intersections are created between $W^u(O)$
and $W^s(O)$,  a discrete Shilnikov attractor is formed containing the fixed point $O$,
Fig.~\ref{shil_or_mirmap}h. In Fig.~\ref{shil_or_homoclinic},  the numerically
obtained stable separatrix $W^{s+}(O)$ is shown which confirms the existence of homoclinic intersections of $W^{u}(O)$ and
$W^{s+}(O)$ within the attractor.
\end{itemize}

Now we consider the case of the appearance of nonorientable Shilnikov attractor again in map \eqref{Miramap}. (Recall that in this map such attractor was found in \cite{KSK21}). Figure~\ref{shil_nonor-y2mir} shows main stages of the development of this attractor when $B = -0.915, A = -2.786$ are fixed and $C$ changes. The stages of the attractor creation follow the general bifurcation scenario described in \cite{KSK21, GGKozS21}, see also section~\ref{sec:Shnor_scen}:
\begin{figure}[ht]
\centerline{\epsfig{file=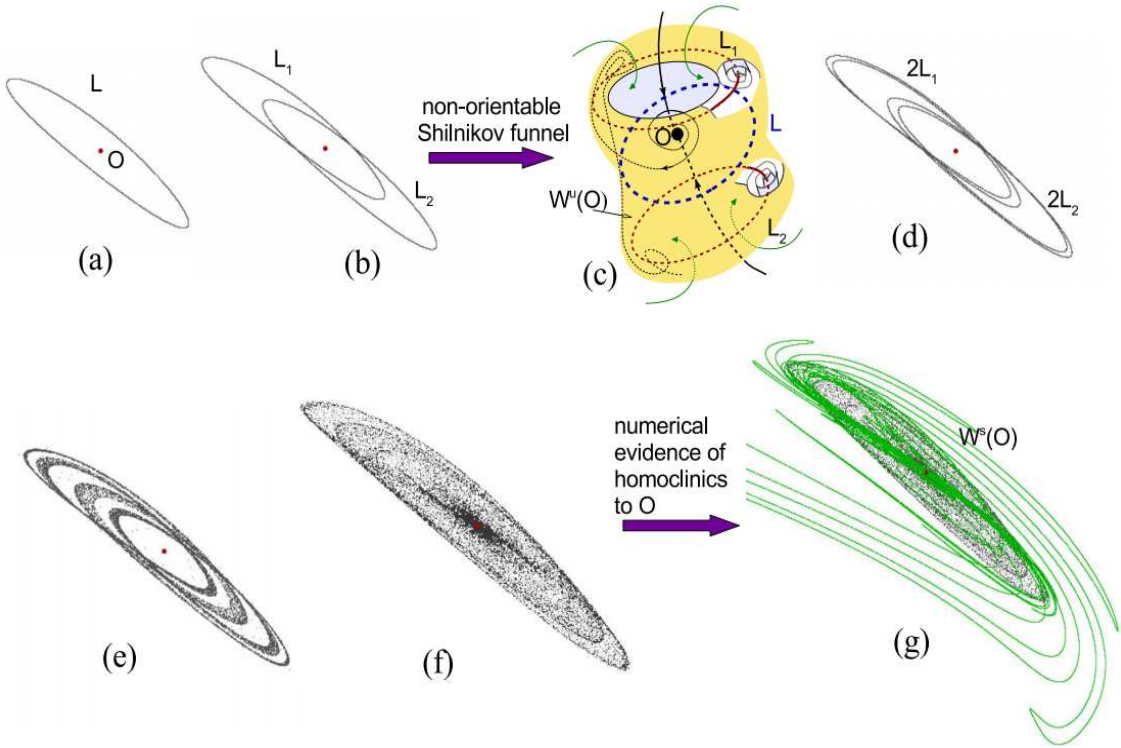, width=16cm }}
\hspace{-1cm}
  \caption{{\footnotesize
Stages of the emergence of a nonorientable Shilnikov attractor in map (\ref{Miramap}) for $B = -0.915, A = -2.786$ when $C$ changes (picture (c) is only schematic). }}
   \label{shil_nonor-y2mir}
\end{figure}
\begin{itemize}
\item
At the beginning, the attractor is the fixed point $O(0,0,0)$ that is stable for
$ 4.701> C> c_1 \approx - 2.712 $).
\item
Then, at $C=c_1$, the point $O$ undergoes a supercritical Andronov-Hopf bifurcation and becomes a nonorientable saddle-focus (1,2) and a stable closed invariant curve $L$ is born (see
Fig.~\ref{shil_nonor-y2mir}a for $ C = - 2.723 $).
\item
Next, the curve $L$ undergoes a component-doubling bifurcation \cite{GGT21}: a pair of stable period-2 curves $L_1$ and $L_2$ are born from $L$ (these curves are shown in Fig.~\ref{shil_nonor-y2mir}b for $ C = -2.733$). With the further change of $C$, the two-dimensional unstable manifold of $O_\mu$ begins to wind up on both the curves $L_1$ and $L_2$ and a nonorientable Shilnikov whirlpool is formed, which has a ``double roll'' shape (it is shown schematically, in Fig.~\ref{shil_nonor-y2mir}c).
\item
Further, the dynamics inside the funnel becomes more complicated. In particular, several bifurcations of doubling of invariant curves occur giving rise of strange attractors inside the whirlpool, (see Fig.~\ref{shil_nonor-y2mir}e, where one of such attractors is shown for $C= -2.736$).
\item
Finally, when homoclinic intersections are created between $W^u(O)$ and $W^s(O)$, a non-orientable discrete Shilnikov attractor is formed, (see Fig.~\ref{shil_nonor-y2mir}f for $C=-2.743$). In Fig.~\ref{shil_nonor-y2mir}g a piece of $W^s(O)$ is shown confirming that the attractor contains the point $O$.
\end{itemize}

\section{Conclusion}

We note that the Shilnikov's ideas about the importance of scenarios for the emergence of spiral attractors and, more generally,   homoclinic attractors (strange attractors containing  either an equilibrium  or a periodic orbit and its entire unstable manifold) were implemented in a number of our works. The first one was our work with L.P. Shilnikov \cite{GGS12}, in which we gave a phenomenological description of scenarios for the emergence of discrete attractors of various types (discrete Lorenz attractors, figure-8 attractors and Shilnikov attractors), and also gave examples of the implementation of these scenarios in one-parameter families of three-dimensional maps. Now the theory of discrete homoclinic attractors looks quite advanced, and it is, of course, richer than the corresponding theory of three-dimensional flows, see e.g. \cite{GGOT13,GGKT14,GG16,GGKozS21,GGKS21}. On the other hand, there is also a certain feedback here: a detailed study of discrete homoclinic attractors also has led to the discovery of new types of attractors of three-dimensional flows (for example, Lorenz attractors with several equilibria \cite{GGKS21}). By the way, these studies are still ongoing.

\section*{Acknowledgements}
This paper was carried out in the framework of the Russian Ministry of Science and Education grant No. 0729-2020-0036.
A. Gonchenko was supported by the RSciF Grant No. 20-71-00079 (Section 4 and 5). Yu. Bakhanova and A. Kazakov were supported by the RSciF Grant No. 19-71-10048. (Section 3).  S. Gonchenko, A. Kazakov and E. Samylina thank the Theoretical Physics and Mathematics Advancement Foundation ``BASIS'', Grant No. 20-7-1-36-5, for support of scientific investigations.


\begin{thebibliography}{99}

\bibitem{Sh65}
Shilnikov L. P. A case of the existence of a denumerable set of periodic motions //Doklady Akademii Nauk. Russian Academy of Sciences, 1965. Vol. 160, no. 3, P. 558--561.
%
\bibitem{Sh62}
Shilnikov L. P. Some cases of generation of periodic motions in n-space //Doklady Akademii Nauk. Russian Academy of Sciences, 1962. Vol. 143, no. 2, P. 289--292.
%
\bibitem{Sh63}
Shilnikov L. P. Some cases of generation of period motions from singular trajectories //Matematicheskii Sbornik. 1963. Vol. 103, no. 4, P. 443-466.
%
\bibitem{Sh68}
Shilnikov L. P. On the generation of a periodic motion from trajectories doubly asymptotic to an equilibrium state of saddle type //Matematicheskii Sbornik. 1968. Vol. 119, no. 3, P. 461--472.
%
\bibitem{Sh70}
Shilnikov L. P. On the question of the structure of an extended neighborhood of a structurally stable state of equilibrium of saddle-focus type //Mat. Sb.(NS). 1970. Vol. 81, no. 123, P. 92--103.
%
\bibitem{Sh84}
Shilnikov L. P. Bifurcation theory and turbulence //Nonlinear and Turbulent Processes in Physics. 1984. P. 1627.
%
\bibitem{Sh86}
Shilnikov L. P. Bifurcation theory and turbulence. I, //Methods Qual. Theory Differ. Equations. 1986. P. 150--163.

\bibitem{ChuaKomMat1986} Chua L.O., Komuro M., Matsumoto T. The double scroll family // Circuits and Systems. IEEE Transactions on. 1986. Vol. 33, no. 11. P. 1072--1118.
%
\bibitem{Ani1990} Anishchenko V.S. Complex oscillations in simple systems // M., 1990 (in Russian).
%
\bibitem{ArrMeuGad1987} Arecchi F.T., Meucci R., Gadomski W. Laser dynamics with competing instabilities // Physical
Review Letters. 1987. Vol. 58, no. 21. P. 2205.
%
\bibitem{ArrLapMeuRov1988} Arecchi F.T., Lapucci A., Meucci R., Roversi J.A., Coullet P.H. Experimental characterization of
Shil'nikov chaos by statistics of return times // EPL (Europhysics Letters). 1988. Vol. 6, no. 8. P. 677.
%
\bibitem{PisurtMeuBrechou3} Pisarchik A.N., Meucci R., Arecchi F.T. Theoretical and experimental study of discrete behavior
of Shilnikov chaos in a CO2 laser // The European Physical Journal D-Atomic, Molecular, Optical
and Plasma Physics. 2001. Vol. 13, no. 3. P. 385--391.
%
\bibitem{ZhouKurtAllBocMeu2003}
 Zhou C.S., Kurths J., Allaria E., Boccaletti S., Meucci R., Arecchi F.T. Constructive effects of
noise in homoclinic chaotic systems // Physical Review E. 2003. Vol. 67, no. 6. P. 066220.
%
\bibitem{ArgArnRich1987}
Argoul F., Arneodo A., Richetti P. Experimental evidence for homoclinic chaos in the BelousovZhabotinskii reaction // Physics Letters A. 1987. Vol. 120, no. 6. P. 269--275.
%
\bibitem{ArgArnEleRich1993} Arneodo A., Argoula F., Elezgarayab J., Richettia P. Homoclinic chaos in chemical systems //
Physica D: Nonlinear Phenomena. 1993. Vol. 62, no. 1. P. 134--169.
%
\bibitem{FeudelPei2000}
Feudel U., Neiman A., Pei X., Wojtenek W., Braun H., Huber M., Moss F. Homoclinic bifurcation
in a Hodgkin-Huxley model of thermally sensitive neurons // Chaos: An Interdisciplinary Journal
of Nonlinear Science. 2000. Vol. 10, no. 1. P. 231--239.
%
\bibitem{PartEdGri2003}
Parthimos D., Edwards D.H., Griffith T.M. Shilnikov homoclinic chaos is intimately related to
type-III intermittency in isolated rabbit arteries: role of nitric oxide // Physical Review E. 2003.
Vol. 67, no. 5. P. 051922.
%
\bibitem{KopGasSlu1992n2001}Koper M.T.M., Gaspard P., Sluyters J.H. Mixed-mode oscillations and incomplete homoclinic
scenarios to a saddle focus in the indium / thiocyanate electrochemical oscillator // Journal of
Chemical Physics. 1992. Vol. 97, no. 11. P. 8250--8260.
%
\bibitem{CheWoaDomn2001}
Chedjou J.C., Woafo P., Domngang S. Shilnikov chaos and dynamics of a self-sustained electromechanical transducer // Journal of vibration and acoustics. 2001. Vol. 123, no. 2. P. 170--174.
%
\bibitem{BasHud1988}
 Bassett M.R., Hudson J.L. Shilnikov chaos during copper electrodissolution // Journal of Physical
Chemistry. 1988. Vol. 92, no. 24. P. 6963--6966.
%
\bibitem{Noh2009}
Noh T. Shilnikov chaos in the oxidation of formic acid with bismuth ion on Pt ring electrode //
Electrochimica Acta. 2009. Vol. 54, no. 13. P. 3657--3661.
%
\bibitem{Rucklidg1993}
Rucklidge A.M. Chaos in a low-order model of magnetoconvection // Physica D: Nonlinear
Phenomena. 1993. Vol. 62, no. 1. P. 323--337.
%
\bibitem{HenLeviOdeh1991}
Henderson M.E., Levi M., Odeh F. The geometry and computation of the dynamics of coupled
pendula // International Journal of Bifurcation and Chaos. 1991. Vol. 1, no. 01. P. 27--50
%
\bibitem{GGNT97}
Gonchenko S. V., Turaev D. V., Gaspard P., Nicolis G. Complexity in the bifurcation structure of homoclinic loops to a saddle-focus //Nonlinearity. 1997. Vol. 10, no. 2, P. 409.
%
\bibitem{ASh83}
Afraimovich V.S., Shilnikov L.P.  Invariant two-dimensional tori,
their breakdown and stochasticity // in Methods of qualitative theory of
differential equations, Gorky, 3-26 (1983) [English translation in
Am. Math. Soc. Transl., Ser. 2, 149, 201--212 (1991).
%
\bibitem{Sh67}
Shilnikov L.P. Existence of a countable set of periodic motions in a four-dimensional space in an extended neighborhood of a saddle-focus //Doklady Akademii Nauk. 1967. Vol. 172, no. 1, P. 54--57.
%
\bibitem{ArnCoulTres1980}
Arneodo A., Coullet P. H., Tresser C. Occurence of strange attractors in three-dimensional Volterra equations //Physics Letters A. 1980. Vol. 79, no. 4, P. 259--263.
%
\bibitem{ArnCoulTres1981}
Arneodo A., Coullet P. H., Tresser C. Possible new strange attractors with spiral structure //Communications in Mathematical Physics. 1981. Vol. 79, P. 573-579.
%
\bibitem{ArnCoulTres1982}
Arneodo, A., P. Coullet, and C. Tresser. Oscillators with chaotic behavior: An illustration of a theorem by Shil'nikov. // Journal of Statistical Physics. 1982. Vol. 27, no. 1, P. 171--182.
%
\bibitem{ArnCoulTres1985}
Arneodo A., Coullet P.H., Spiegel E.A., Tresser C. Asymptotic chaos. //Physica D. 1985. Vol. 14, P. 327--347.
%
\bibitem{GaspardNicolis83}
Gaspard P., Nicolis G. What can we learn from homoclinic orbits in chaotic dynamics? //Journal of statistical physics. 1983. Vol. 31, no. 3, P. 499--518.
%
\bibitem{LetDutMah1995}
Letellier, C., Dutertre, P. and Maheu, B. Unstable periodic orbits and templates of the R\"ossler system: toward a systematic topological characterization. Chaos: An Interdisciplinary Journal of Nonlinear Science. 1995 Vol. 5, no. 1, P. 271-282.
%
\bibitem{Malykh20}
Malykh, S., Bakhanova, Y., Kazakov, A., Pusuluri, K., Shilnikov, A. Homoclinic chaos in the R\"ossler model. Chaos: An Interdisciplinary Journal of Nonlinear Science, 2020 Vol. 30, no. 11, p. 113126.
%
\bibitem{Gallas2010}
Gallas J. A. The structure of infinite periodic and chaotic hub cascades in phase diagrams of simple autonomous flows //Int. J. Bifurcation Chaos. 2010. Vol. 20, no. 02, P. 197--211.
%
\bibitem{GGT21}
Gonchenko A. S., Gonchenko S. V., Turaev D. V. Doubling of invariant curves and chaos in three-dimensional diffeomorphisms //Chaos. 2021. Vol. 31, P. 113130
\bibitem{GGS12}
Gonchenko A.S., Gonchenko S.V., Shilnikov L.P. Towards scenarios of chaos appearance in three-dimensional maps// Rus. J. Nonlinear Dynamics. 2012. Vol. 8, P. 3-28.
%
\bibitem{GGKT14}
Gonchenko A. S., Gonchenko, S. V., Kazakov, A. O., Turaev, D. V. Simple scenarios of onset of chaos in three-dimensional maps // International Journal of Bifurcation and Chaos. 2014. Vol. 24, no. 08, P. 1440005

\bibitem{KSK21}
Karatetskaia E., Shykhmamedov A., Kazakov A. Shilnikov attractors in three-dimensional orientation-reversing maps // Chaos: An Interdisciplinary Journal of Nonlinear Science. 2021. Vol. 31, no. 1, P. 011102.

\bibitem{GGKozS21}
Gonchenko A. S., Gonchenko M. S., Kozlov A. D., Samylina E.A. On scenarios of the onset of homoclinic attractors in three-dimensional non-orientable maps // Chaos: An Interdisciplinary Journal of Nonlinear Science. 2021. Vol. 31, no. 4, P. 043122.

\bibitem{GGOT13}
Gonchenko S. V., Gonchenko A. S., Ovsyannikov I. I., Turaev D. Examples of Lorenz-like attractors in Henon-like maps //Mathematical Modelling of Natural Phenomena. 2013. Vol. 8, no. 5, P. 48--70.

\bibitem{GG16}
Gonchenko A. S., Gonchenko S. V. Variety of strange pseudohyperbolic attractors in threedimensional generalized Henon maps // Physica D. 2016. Vol. 337,  P. 43--57.

\bibitem{GGKS21}
Gonchenko S. V., Gonchenko A. S., Kazakov, A. O., Samylina E. A. On discrete Lorenz-like attractors //Chaos: An Interdisciplinary Journal of Nonlinear Science. 2021. Vol. 31, no. 2, P. 023117.

\end{thebibliography}
\end{document}